\documentclass[12pt]{amsart}

\DeclareRobustCommand{\eulerian}{\genfrac<>{0pt}{}}
\usepackage{graphicx}
\usepackage{tikz}
\usetikzlibrary{arrows,shapes,automata,snakes,backgrounds,petri,positioning}
\usepackage[utf8]{inputenc}
\usepackage{mathtools}

\newtheorem{thm}{Theorem}
\newtheorem{cl}{Claim}

\newtheorem{cor}{Corollary}

\theoremstyle{definition}

\DeclareMathOperator{\maj}{maj}
\DeclareMathOperator{\des}{des}

\title{A Stochastic Approach to Eulerian Numbers}

\author{Kiana Mittelstaedt}

\begin{document}

\maketitle

\begin{abstract}
We examine the aggregate behavior of one-dimensional random walks in a model known as (one-dimensional) Internal Diffusion Limited Aggregation. In this model, a sequence of $n$ particles perform random walks on the integers, beginning at the origin. Each particle walks until it reaches an unoccupied site, at which point it occupies that site and the next particle begins its walk. After all walks are complete, the set of occupied sites is an interval of length $n$ containing the origin. We show the probability that $k$ of the occupied sites are positive is given by an Eulerian probability distribution. Having made this connection, we use generating function techniques to compute the expected run time of the model.
\end{abstract}

\section*{Introduction}

\begin{figure}
\[
\begin{tikzpicture}
\tikzstyle{wn}=[inner sep =3, fill=white, draw=black, circle];
\tikzstyle{bn}=[inner sep =3, fill=black, draw=black, circle];
\draw (0,0)--(1,0) node[wn] {} --(2,0) node[wn] {}--(3,0) node[wn] {}--(4,0) node[bn] {}--(5,0) node[wn] {}--(6,0) node[wn] {}--(7,0);
\draw (4,.2) node[above] {$\circ$}; \draw (4,.45) node[above] {$\circ$}; \draw (4,.7) node[above] {$\circ$}; 
\draw (4,1) node[above] {$\vdots$};
\draw (4,-.2) node[below] {$\star$};
\draw (0,0) node[left] {$\cdots$};
\draw (7,0) node[right] {$\cdots$};
\end{tikzpicture}
\]
\caption{The starting graph for our coin-tossing game.}\label{fig:begstate}
\end{figure}
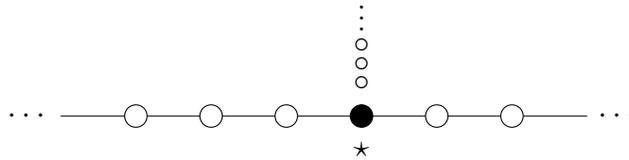

Consider the following coin-tossing game. To play the game, we will begin with a graph that has one active node with infinitely many inactive nodes on either side of it. In this paper, we will also refer to active nodes as occupied sites and inactive nodes as unoccupied sites. We will label the initial occupied site the origin. At the origin, there is a queue of particles that will perform random walks on the graph. This is seen in Figure \ref{fig:begstate}, where occupied sites are black, unoccupied sites are white, the origin is labeled with ``$\star$'' and the particles in the queue are indicated with ``$\circ$''.

We imagine dropping the first particle in the queue onto the origin. Flipping heads will move the particle one position to the right and flipping tails will move it one position to the left. The game consists of flipping a coin repeatedly until the particle reaches an unoccupied site. Once an unoccupied site is reached, the particle settles there and occupies it. Now, the next particle in the queue is dropped onto the origin and the game begins again.

\begin{figure}
\[
\begin{tikzpicture}
\draw (0,0) node[scale=.75] (a) {
\begin{tikzpicture}
\draw (2.5,0) -- (3,0) node[inner sep =3, fill=white, draw=black, circle] {}--(4,0) node[inner sep =3, fill=black, draw=black, circle] {}--(5,0) node[inner sep =3, fill=white, draw=black, circle] {} -- (5.5,0);
\draw (4,-.2) node[below] {$\star$}; \draw (4,.2) node[above] {$\circ$}; \draw (2.5,0) node[left] {$\dots$}; \draw (5.5,0) node[right] {$\dots$};

\end{tikzpicture}
};
\draw (0,-2) node[scale=.75] (a2) {
\begin{tikzpicture}
\draw (2.5,0) -- (3,0) node[inner sep =3, fill=white, draw=black, circle] {}--(4,0) node[inner sep =3, fill=black, draw=black, circle] {}--(5,0) node[inner sep =3, fill=white, draw=black, circle] {} -- (5.5,0);
\draw (4,-.2) node[below] {$\star$}; \draw (5,.2) node[above] {$\circ$}; \draw (2.5,0) node[left] {$\dots$}; \draw (5.5,0) node[right] {$\dots$};
\end{tikzpicture}
};
\draw (5.5,-2) node[scale=.75] (b) {
\begin{tikzpicture}
\draw (2.5,0) -- (3,0) node[inner sep =3, fill=white, draw=black, circle] {}--(4,0) node[inner sep =3, fill=black, draw=black, circle] {}--(5,0) node[inner sep =3, fill=black, draw=black, circle] {}--(6,0) node[inner sep =3, fill=white, draw=black, circle] {} -- (6.5,0);
\draw (4,-.2) node[below] {$\star$}; \draw (4,.2) node[above] {$\circ$}; \draw (2.5,0) node[left] {$\dots$}; \draw (6.5,0) node[right] {$\cdots$};
\end{tikzpicture}
};
\draw (5.5,-4) node[scale=.75] (b2) {
\begin{tikzpicture}
\draw (2.5,0) -- (3,0) node[inner sep =3, fill=white, draw=black, circle] {}--(4,0) node[inner sep =3, fill=black, draw=black, circle] {}--(5,0) node[inner sep =3, fill=black, draw=black, circle] {}--(6,0) node[inner sep =3, fill=white, draw=black, circle] {} -- (6.5,0);
\draw (4,-.2) node[below] {$\star$}; \draw (5,.2) node[above] {$\circ$}; \draw (2.5,0) node[left] {$\dots$}; \draw (6.5,0) node[right] {$\cdots$};
\end{tikzpicture}
};
\draw (5.5,-6) node[scale=.75] (b3) {
\begin{tikzpicture}
\draw (2.5,0) -- (3,0) node[inner sep =3, fill=white, draw=black, circle] {}--(4,0) node[inner sep =3, fill=black, draw=black, circle] {}--(5,0) node[inner sep =3, fill=black, draw=black, circle] {}--(6,0) node[inner sep =3, fill=white, draw=black, circle] {} -- (6.5,0);
\draw (4,-.2) node[below] {$\star$}; \draw (4,.2) node[above] {$\circ$}; \draw (2.5,0) node[left] {$\dots$}; \draw (6.5,0) node[right] {$\cdots$};
\end{tikzpicture}
};
\draw (5.5,-8) node[scale=.75] (b4) {
\begin{tikzpicture}
\draw (2.5,0) -- (3,0) node[inner sep =3, fill=white, draw=black, circle] {}--(4,0) node[inner sep =3, fill=black, draw=black, circle] {}--(5,0) node[inner sep =3, fill=black, draw=black, circle] {}--(6,0) node[inner sep =3, fill=white, draw=black, circle] {} -- (6.5,0);
\draw (4,-.2) node[below] {$\star$}; \draw (3,.2) node[above] {$\circ$}; \draw (2.5,0) node[left] {$\dots$}; \draw (6.5,0) node[right] {$\cdots$};
\end{tikzpicture}
};
\draw (11.5,-8) node[scale=.75] (c) {
\begin{tikzpicture}
\draw (1.5,0) -- (2,0) node[inner sep =3, fill=white, draw=black, circle] {}--(3,0) node[inner sep =3, fill=black, draw=black, circle] {}--(4,0) node[inner sep =3, fill=black, draw=black, circle] {}--(5,0) node[inner sep =3, fill=black, draw=black, circle] {}--(6,0) node[inner sep =3, fill=white, draw=black, circle] {} -- (6.5,0);
\draw (4,-.2) node[below] {$\star$}; \draw (4,.2) node[above] {$\circ$}; \draw (1.5,0) node[left] {$\dots$}; \draw (6.5,0) node[right] {$\cdots$};
\end{tikzpicture}
};
\draw[->] (a) -- node[midway,right] {H} (a2);
\draw[->] (b) -- node[midway,right] {H} (b2);
\draw[->] (b2) -- node[midway,right] {T} (b3);
\draw[->] (b3) -- node[midway,right] {T} (b4);
\draw[->] (a2) -- node[midway,above,sloped] {new particle} (b);
\draw[->] (b4) -- node[midway,above,sloped] {new particle} (c);
\end{tikzpicture}
\] 
\caption{Playing the coin-tossing game.}\label{fig:game}
\end{figure}
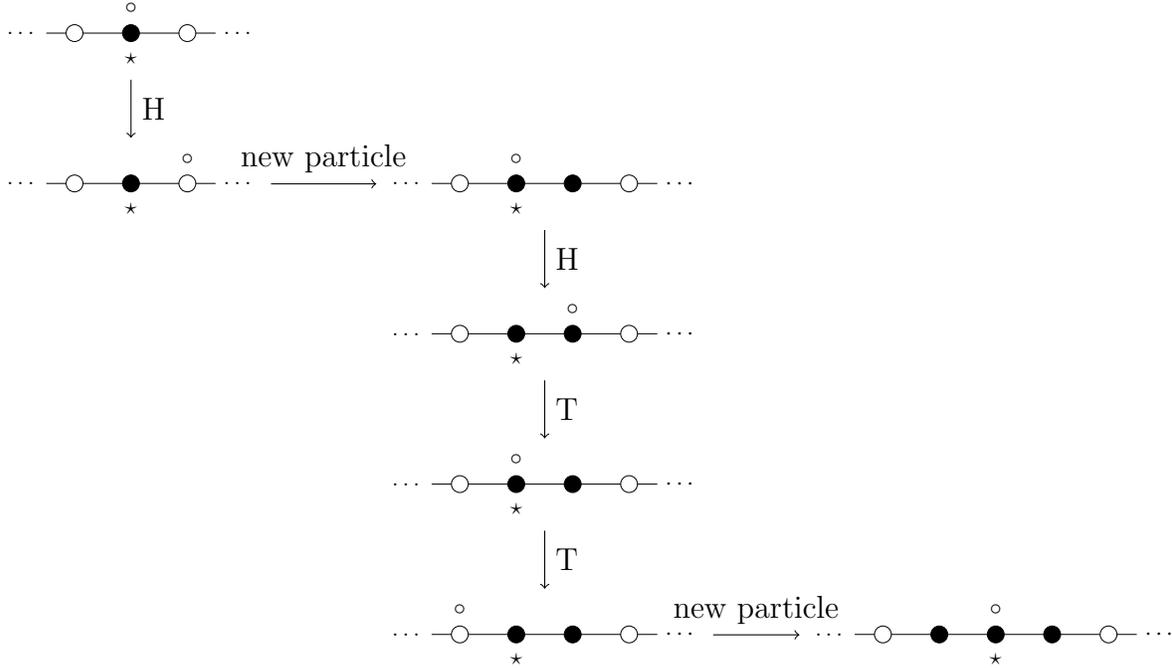

The illustration in Figure \ref{fig:game} shows an example of the game being played for the first two particles in the queue. 
In the example, we first flip heads. Thus, the first particle moves to the unoccupied site immediately to the right of the origin, settles there, and occupies it. Now, the graph has two occupied sites, each with infinitely many unoccupied sites on each side. Since a new site was occupied, the second particle in the queue is dropped onto the origin and the game begins again. With the next toss, there are two unoccupied sites that could become occupied: the site to the left of the origin or the site two to the right of the origin. Suppose the next coin toss is heads, followed by tails, then tails again. These three tosses would move the new particle to the occupied site, back to the origin, and then to the unoccupied site to the left of the origin. As a result, the particle settles at this site and occupies it. Next, the third next particle in the queue is dropped on the origin and the process continues.

In our game, each particle on its own is performing a simple one-dimensional random walk, which is a well-studied stochastic model. See \cite{Sauer}. This game we have described, of sending particles to take random walks in succession, is known as Internal Diffusion Limited Aggregation (Internal DLA). This process was defined by Meakin and Deutch in \cite{MeakinDeutch} to serve as a model for ``chemical etching.'' In a more mathematical work by Diaconis and Fulton \cite{DiaconisFulton}, Internal DLA was studied using the language of Markov chains, with interesting algebraic connections. The limiting behavior of this model has been frequently studied in dimensions two and three (e.g., in a square grid, the limiting shape is a ball), but here we will consider the one-dimensional case.     

Some potential questions we can ask in the one-dimensional case for Internal DLA are:
\begin{enumerate}
  \item What is the probability that, after releasing $n$ particles, $k$ of them settle to the right of the origin? 
  \item On average, how many coin tosses will it take for $n$ particles to settle? 
  \item Given a large queue of particles, if the coin is tossed $N$ times, what is the probability that there are $n$ occupied sites? 
  \item What if the particles do not perform a simple random walk and instead there is a probability assigned to each site? 
\end{enumerate}
In this paper, we prove results for questions $(1)$ and $(2)$. Diaconis and Fulton addressed question $(1)$ and found the limiting distribution (as $n \to \infty$) is normal \cite[Proposition 3.2]{DiaconisFulton}. We give a more precise characterization here. While we have investigated questions $(3)$ and $(4)$, we do not yet have any major results. Some comments on further pursuit of these questions are included at the end of the paper. 

In Section 1, we consider question $(1)$ as a probability distribution for our coin-tossing game. That is, given that $n$ particles have left the origin and settled, we determine the probability that $k$ are to the right of the origin. We first used a Monte Carlo simulation to predict a distribution. Each simulation ran $100\,000$ trials of the game and returned the minimum and maximum positions reached on the graph (assuming the origin is at position 0). Table \ref{tab:Monte1} shows the observed probability of the game reaching maximum ending position $k$ for $1 \leq n \leq 7$.

\begin{table}
\[
\begin{array}{c | c c c c c c c | c}
n \backslash k & 0 & 1 & 2 & 3 & 4 & 5 & 6 & SSE\\
\hline
1 & 1.00000 &&&&&&& 0.00000 \\
2 & 0.50266 & 0.49734 &&&&&& 1.41512\mathrm{e}{-10} \\
3 & 0.16730 & 0.66561 & 0.16709 &&&&& 1.69687\mathrm{e}{-10} \\
4 & 0.04235 & 0.45694 & 0.45763 & 0.04308 &&&& 4.90051\mathrm{e}{-10}\\
5 & 0.00833 & 0.21704 & 0.54711 & 0.21906 & 0.00846 &&& 1.42356\mathrm{e}{-10}\\
6 & 0.00131 & 0.08109 & 0.41734 & 0.41990 & 0.07903 & 0.00133 && 8.36380\mathrm{e}{-10}\\
7 & 0.00018 & 0.02418 & 0.23675 & 0.47821 & 0.23671 & 0.02381 & 0.00016 & 1.82768\mathrm{e}{-10}
\end{array}
\]
\caption{Monte Carlo results for reaching ending position $k$ for $1 \leq n \leq 7$. }\label{tab:Monte1}
\end{table}

We know that when $n=1, k=0$ with probability $1$ because the only occupied site is the origin. Furthermore, we know that when $n=2,$ both $k=0$ and $k=1$ occur with probability $\frac{1}{2}$ because $50\%$ of the time we will flip heads and settle at the site to the right of the origin, and $50\%$ of the time we will flip tails and settle at the site to the left of the origin. The data in Table \ref{tab:Monte1} supports this fact. 

When $n=3$, our data suggests that there is a common denominator of $6$, and when $n=4$, our data suggests the common denominator may be $24$. For example, when $n=3$ and $k=0$, we see a probability of $0.16730$ which is very close to $0.1\overline{666} = \frac{1}{6}$. Similarly, when $n=4$ and $k=0$, we see $0.04235$ which is very close to $0.041\overline{6} = \frac{1}{24}$. If we consider these observations with our known probabilities for $n=1$ and $n=2$, we see that our denominators for $1 \leq n \leq 4$ may be $n!$. 

Upon recognizing this, we multiplied each row of Table \ref{tab:Monte1} by a common denominator of $n!$ to try and find a nice distribution for the numerators. For example, multiplying row $5$ by $5!$ yields
\[
[0.9996 \quad 26.0448 \quad 65.6532 \quad 26.2872 \quad 1.0152].
\]
After multiplication by $n!$, we found that our data looked very close to the triangle of Eulerian numbers. See Table \ref{table:eulnum}. Upon comparing our data in Table \ref{tab:Monte1} to the corresponding Eulerian number divided by $n!$, we got a very close fit. The sum of squared errors, denoted $SSE$, for each $1 \leq n \leq 7$ is shown in the last column in Table \ref{tab:Monte1}, comparing our data to the corresponding Eulerian number divided by $n!$. Thus, our Monte Carlo simulation suggests that our coin-tossing game follows an Eulerian probability distribution and we will prove this result in Theorem \ref{thm:eulprob}.\footnote{As communicated by Lionel Levine, it seems Jim Propp made this connection in unpublished work. In fact it was Levine, in a visit to DePaul University, who suggested we investigate this question in the first place.}

Recall that for any $0 \leq k \leq n-1,$ the \textit{Eulerian number $\eulerian{n}{k}$} is the number of permutations of $\{1,2,\dots,n\}$ with $k$ descents. See \cite{Petersen}. The Eulerian numbers satisfy the recurrence
\begin{equation}\label{eq:eulrec}
\eulerian{n}{k} = (n - k)\eulerian{n - 1}{k - 1} + (k + 1)\eulerian{n - 1}{k},
\end{equation}
for any $n \geq 1$ and $0 \leq k \leq n - 1$, with boundary values $\eulerian{n}{0} = \eulerian{n}{n-1} = 1$ for all $n\geq 0$. Table \ref{table:eulnum} shows the triangle of Eulerian numbers for $1 \leq n \leq 6$. The Eulerian numbers are a finite combinatorial distribution and our coin-tossing game is a stochastic process. How these mathematical concepts are related is very puzzling at first glance. 
\begin{table}
\begin{tabular}{c | c c c c c c}
$n \backslash k$ & 0 & 1 & 2 & 3 & 4 & 5\\
\hline
1 & 1 \\
2 & 1 & 1 \\
3 & 1 & 4 & 1 \\
4 & 1 & 11 & 11 & 1 \\
5 & 1 & 26 & 66 & 26 & 1 \\
6 & 1 & 57 & 302 & 302 & 57 & 1
\end{tabular}
\caption{Triangle of the Eulerian numbers for $1 \leq n \leq 6$.}
\label{table:eulnum}
\end{table}

In Section 2, we consider question $(2)$, the expected run time of the game. From our Monte Carlo simulations we recorded the average number of coin tosses needed for $n$ particles to settle. Let $E_n$ denote this number of expected tosses. The data from our simulation is shown in Table \ref{tab:Monte2}. 

\begin{table}
\[
\begin{array}{c | c ccccccccccccc}
n& 1 & 2 & 3 & 4& 5& 6& 7& 8& 9 & 10  & 15 & 20 \\
\hline
E_n & 0 & 1 & 3.00 & 6.66 & 12.49 & 20.98 & 32.64 & 47.96 & 67.54 & 91.68 & 300.24 & 699.79
\end{array}
\]
\caption{Monte Carlo results for $E_n$ for $1 \leq n \leq 10$, $n=15,$ and $n=20$.}\label{tab:Monte2}
\end{table}

Using Lagrange interpolation, we find the number of tosses needed for the game to have $n$ occupied sites is roughly
\[
\frac{1}{12} n^3 + \frac{1}{12} n^2.
\]
We will prove this result for expected run time in Theorem \ref{thm:tosstosn}. 

\section{Probability Distribution}

To begin, we define the possible states of the game and a two-parameter family of probabilities relating these states. 

For any $n \geq 1$ and $0 \leq k \leq n - 1$, let $s(n,k)$ denote the graph with $n$ occupied sites, $k$ of which are to the right of the origin. State $s(n,k)$ is shown in Figure \ref{fig:snk}. We think of these graphs as possible states of our coin-tossing game. For example, 
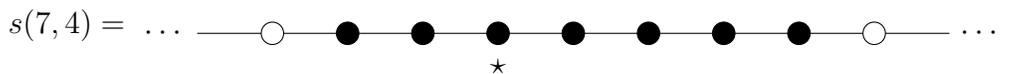
\[
s(7,4) = \begin{tikzpicture}[baseline=0]
\tikzstyle{wn}=[inner sep =3, fill=white, draw=black, circle];
\tikzstyle{bn}=[inner sep =3, fill=black, draw=black, circle];
\draw (0,0)--(1,0) node[wn] {} -- (2,0) node[bn] {}--(3,0) node[bn] {}--(4,0) node[bn] {}--(5,0) node[bn] {}--(6,0) node[bn] {}--(7,0) node[bn] {}--(8,0) node[bn] {}--(9,0) node[wn] {}--(10,0);
\draw (4,-.2) node[below] {$\star$}; \draw (0,0) node[left] {$\cdots$}; \draw (10,0) node[right] {$\cdots$};
\end{tikzpicture}.
\]

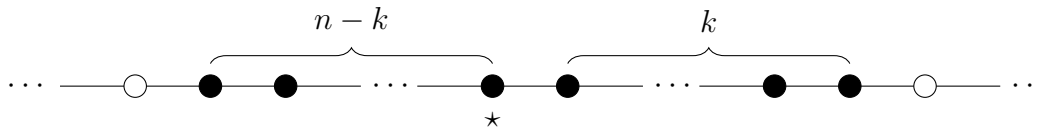
\begin{figure}
\[
\begin{tikzpicture}[baseline=0]
\tikzstyle{wn}=[inner sep =3, fill=white, draw=black, circle];
\tikzstyle{bn}=[inner sep =3, fill=black, draw=black, circle];
\draw (0,0)--(1,0) node[wn] {} -- (2,0) node[bn] {}--(3,0) node[bn] {}--(4,0);
\draw (4.75,0)--(5.75,0) node[bn] {}--(6.75,0) node[bn] {}--(7.75,0);
\draw (8.5,0)--(9.5,0) node[bn] {}--(10.5,0) node[bn] {}--(11.5,0) node[wn] {}-- (12.5,0);
\draw (0,0) node[left] {$\cdots$}; \draw (4,0) node[right] {$\cdots$}; \draw (7.75,0) node[right] {$\cdots$}; \draw (12.5,0) node[right] {$\cdots$};
\draw (5.75,-.2) node[below] {$\star$};
\draw [decoration={brace,amplitude=6pt},decorate] (6.75,0.3) -- (10.5,0.3) node[midway,above=8pt] {$k$};
\draw [decoration={brace,amplitude=6pt},decorate] (2,0.3) -- (5.75,0.3) node[midway,above=8pt] {$n-k$};
\end{tikzpicture}
\]
\caption{State $s(n,k)$.}\label{fig:snk}
\end{figure}

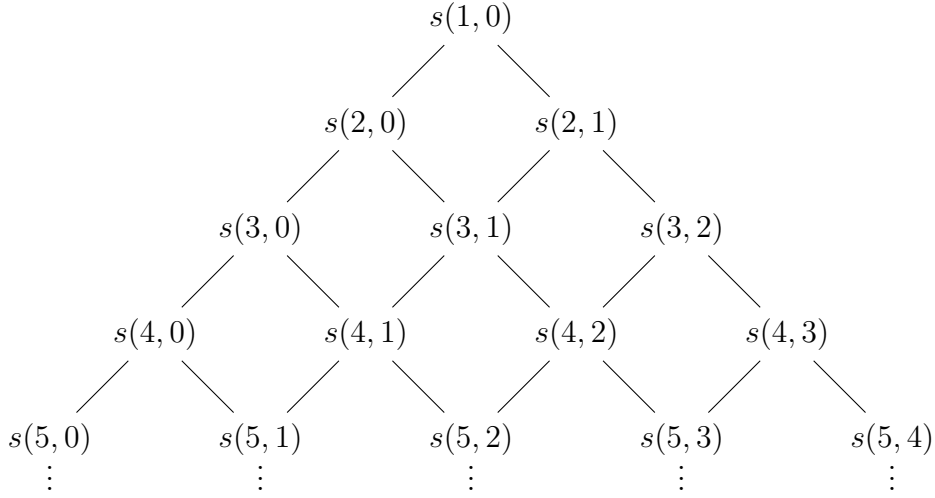
\begin{figure}
\begin{center}
\[
\begin{tikzpicture}[scale=.7]
  \node (oz) at (0,4) {$s(1,0)$}; 
  \node (tz) at (-2,2) {$s(2,0)$}; \node (two) at (2,2) {$s(2,1)$};
  \node (thz) at (-4,0) {$s(3,0)$}; \node (tho) at (0,0) {$s(3,1)$}; \node (thtw) at (4,0) {$s(3,2)$};
  \node (fz) at (-6,-2) {$s(4,0)$}; \node (fo) at (-2,-2) {$s(4,1)$}; \node (ftw) at (2,-2) {$s(4,2)$}; \node (fth) at (6,-2) {$s(4,3)$};
  \node (ffz) at (-8,-4) {$s(5,0)$}; \node (ffo) at (-4,-4) {$s(5,1)$}; \node (fftw) at (0,-4) {$s(5,2)$}; \node (ffth) at (4,-4) {$s(5,3)$}; \node (fff) at (8,-4) {$s(5,4)$};
  \draw (oz) -- (tz); \draw (oz) -- (two); 
  \draw (tz) -- (thz); \draw (tz) -- (tho); \draw (two) -- (tho); \draw (two) -- (thtw);
  \draw (thz) -- (fz); \draw (thz) -- (fo); \draw (tho) -- (fo); \draw (tho) -- (ftw); \draw (thtw) -- (ftw); \draw (thtw) -- (fth);
  \draw (fz) -- (ffz); \draw (fz) -- (ffo); \draw (fo) -- (ffo); \draw (fo) -- (fftw); \draw (ftw) -- (fftw); \draw (ftw) -- (ffth); \draw (fth) -- (ffth); \draw (fth) -- (fff);
  \draw (ffz) node[below] {\vdots}; \draw (ffo) node[below] {\vdots}; \draw (fftw) node[below] {\vdots}; \draw (ffth) node[below] {\vdots}; \draw (fff) node[below] {\vdots};
\end{tikzpicture}
\]
\end{center}
\caption{Possible states of the game for $0 \leq n \leq 5$.}\label{fig:states}
\end{figure}

The graph $s(n,k)$ can only be reached from $s(n - 1,k - 1)$ or $s(n - 1,k)$ for any $n \geq 1$ and $1 \leq k \leq n-2$. For example, state $s(4,2)$ can be reached only from state $s(3,1)$ or state $s(3,2)$. The boundary cases, states $s(n,0)$ and $s(n, n - 1)$, are exceptions: $s(n,0)$ can only be reached from $s(n-1,0)$ and $s(n,n-1)$ can only be reached from $s(n-1,n-2)$. The possible states of the game for $0 \leq n \leq 5$ are shown in Figure \ref{fig:states}, with edges to indicate which states can be obtained from another.

While we know that the game will transition from one state to the next by moving down the lattice in Figure \ref{fig:states}, we would like to know how likely each transition is. For example, we can ask how likely it is that the game will go from state $s(2,1)$ to state $s(3,1)$.

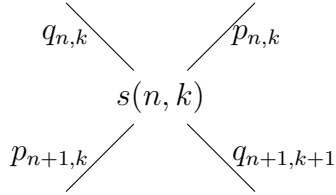
\begin{figure}
\[
\begin{tikzpicture}[scale=.7]
 \node (nk) at (2,0) {$s(n,k)$};
 \node (qnk) at (0,2) {}; 
 \node (pnk) at (4,2) {};
 \node (pn1k) at (0,-2) {};
 \node (qn1k1) at (4,-2) {};
 \draw (qnk) -- node[left]{$q_{n,k}$} (nk); \draw (pnk) -- node[right]{$p_{n,k}$} (nk); \draw (pn1k) -- node[left]{$p_{n+1,k}$}(nk); \draw (qn1k1) -- node[right]{$q_{n+1,k+1}$}(nk);
 \end{tikzpicture}
 \]
 \caption{Transition probabilities for state $s(n,k)$.}\label{fig:snktrans}
 \end{figure}

For any $n \geq 1$ and $0 \leq k \leq n - 1$, let $p_{n,k}$ denote the probability of transitioning from $s(n - 1,k)$ to $s(n,k)$. Similarly, let $q_{n,k}$ denote the probability of transitioning from $s(n - 1,k - 1)$ to $s(n,k)$. We will call $p_{n,k}$ and $q_{n,k}$ transition probabilities. The transition probabilities in and out of state $s(n,k)$ are shown in Figure \ref{fig:snktrans}. 

To help us calculate the transition probabilities, we consider the classic Gambler's Ruin problem. 

\subsection{Gambler's Ruin}\label{gambruin}

In the Gambler's Ruin problem, we consider two gamblers: Gambler A and Gambler B. The gamblers are playing a game where they toss a coin. If the coin flip is tails, Gambler A takes a dollar from Gambler B. Similarly, if heads, Gambler B takes a dollar from Gambler A. Each player starts with some money, say $k$ dollars for Gambler A and $l$ dollars for Gambler B. The game ends when either Gambler A wins all $k+l$ dollars (Gambler B is ruined), or when Gambler B wins all $k+l$ dollars (Gambler A is ruined). The result says that the probability that Gambler A wins is $\frac{k}{k+l}$ and the probability Gambler B wins is $\frac{l}{k+l}$. See \cite{GrinsteadSnell}.

This relates to our model in a straightforward way. Let's assume Gambler A has $4$ dollars and Gambler B has $3$ dollars. This corresponds to the following graph: 

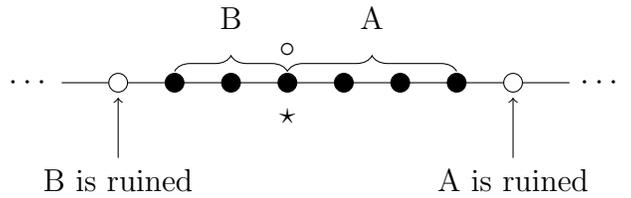
\[
\begin{tikzpicture}[baseline=0]
\tikzstyle{wn}=[inner sep =2.5, fill=white, draw=black, circle];
\tikzstyle{bn}=[inner sep =2.5, fill=black, draw=black, circle];
\draw (0,0)--(.75,0) node[wn] {} -- (1.5,0) node[bn] {}--(2.25,0) node[bn] {}--(3,0) node[bn] {}--(3.75,0) node[bn] {}--(4.5,0) node[bn] {}--(5.25,0) node[bn] {}--(6,0) node[wn] {}--(6.75,0);
\draw (3,-.2) node[below] {$\star$}; \draw (3,.2) node[above] {$\circ$}; \draw (0,0) node[left] {$\cdots$}; \draw (6.75,0) node[right] {$\cdots$};
\draw [->] (.75,-1) node[below] {B is ruined} -- (.75,-.2); \draw [->] (6,-1) node[below] {A is ruined} -- (6,-.2);
\draw [decoration={brace,amplitude=6pt},decorate] (1.5,0.15) -- (3,0.15) node[midway,above=12pt] {B};
\draw [decoration={brace,amplitude=6pt},decorate] (3,0.15) -- (5.25,0.15) node[midway,above=12pt] {A};
\end{tikzpicture}.
\]
We see that the location of the particle on the graph keeps track of the fortune for each player. That is, the money Gambler A has is equal to the number of active nodes to the right of and including the particle, while Gambler B's fortune equals the number of active nodes to the left of and including the particle. 

If the coin toss is heads, Gambler B will take a dollar from Gambler A. So, continuing with our example, Gambler B will have $4$ dollars and Gambler A will have $3$ dollars. Thus, the game moves to the following state:
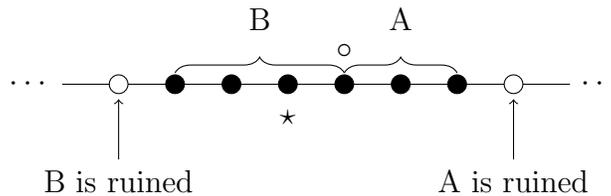
\[
\begin{tikzpicture}[baseline=0]
\tikzstyle{wn}=[inner sep =2.5, fill=white, draw=black, circle];
\tikzstyle{bn}=[inner sep =2.5, fill=black, draw=black, circle];
\draw (0,0)--(.75,0) node[wn] {} -- (1.5,0) node[bn] {}--(2.25,0) node[bn] {}--(3,0) node[bn] {}--(3.75,0) node[bn] {}--(4.5,0) node[bn] {}--(5.25,0) node[bn] {}--(6,0) node[wn] {}--(6.75,0);
\draw (3,-.2) node[below] {$\star$}; \draw (3.75,.2) node[above] {$\circ$}; \draw (0,0) node[left] {$\cdots$}; \draw (6.75,0) node[right] {$\cdots$};
\draw [->] (.75,-1) node[below] {B is ruined} -- (.75,-.2); \draw [->] (6,-1) node[below] {A is ruined} -- (6,-.2);
\draw [decoration={brace,amplitude=6pt},decorate] (1.5,0.15) -- (3.75,0.15) node[midway,above=12pt] {B};
\draw [decoration={brace,amplitude=6pt},decorate] (3.75,0.15) -- (5.25,0.15) node[midway,above=12pt] {A};
\end{tikzpicture}
\]

In general, state $s(n-1,k-1)$ corresponds to Gambler A having $k$ dollars and Gambler B having $l$ dollars such that $k+l=n$. Thus we see that Gambler A wins if the game goes to state $s(n,k-1)$ and Gambler B wins if the game goes to state $s(n,k)$. The Gambler's Ruin problem tells us that Gambler A will win with probability $\frac{k}{n}$ and Gambler B will win with probability $\frac{n-k}{n}$. Thus, we conclude for any $n \geq 1$ and $1 \leq k \leq n - 1$, 
\begin{equation}\label{eq:pqtransprobs}
p_{n,k-1} = \frac{k}{n} \quad \mbox{ and } \quad q_{n,k} = \frac{n-k}{n}.
\end{equation} 

 \subsection{Eulerian distribution}

We use these equations for our transition probabilities to obtain our first main result, which is the probability of starting the game at $s(1,0)$ and arriving at state $s(n,k)$. Let $P(n,k)$ denote this probability. For example, $P(1,0) = 1,$ and since we move immediately to $s(2,0)$ or $s(2,1)$ after our first coin toss, we see $P(2,0) = \frac{1}{2} = P(2,1)$. For $n \geq 3,$ things are more interesting.  

First, we will consider the transition to the boundary states (when $k=0$ or $k=n-1$). We know that state $s(n,0)$ can only be reached from $s(n-1,0)$, so $P(n,0)$ will only depend on $P(n-1,0)$. Similarly, state $s(n,n-1)$ can only be reached from $s(n-1,n-2)$, so $P(n,n-1)$ depends only on $P(n-1,n-2)$. Using our formulas for the transition probabilities in Equation \eqref{eq:pqtransprobs}, we have 
\[
P(n,0) = p_{n,0} P(n-1,0) \quad \mbox{and} \quad P(n,n-1) = q_{n,n-1} P(n-1,n-2).
\]
Since $p_{n,0}=q_{n,n-1}= \frac{1}{n}$ and $P(1,0)=1$, we find 
\begin{equation}\label{eq:boundary}
P(n,0) = \frac{P(n-1,0)}{n} = \frac{1}{n!} \quad \mbox{ and } \quad P(n,n-1) = \frac{P(n-1,n-2)}{n} = \frac{1}{n!}.
\end{equation}

Next, we will look at the non-boundary cases. Because $s(n,k)$ can only be reached from states $s(n-1,k-1)$ and $s(n-1,k)$, we know that $P(n,k)$ will only depend on $P(n-1,k-1)$ and $P(n-1,k)$. Using our formulas for the transition probabilities in Equation \eqref{eq:pqtransprobs}, we find for $1 \leq k \leq n-2$, the following recurrence for $P(n,k)$:
\begin{align}
P(n,k) &= q_{n,k} P(n-1,k-1) + p_{n,k} P(n-1,k), \nonumber \\
&= \frac{n - k}{n} \cdot P(n - 1,k - 1) + \frac{k + 1}{n} \cdot P(n - 1,k) \label{eq:recursiveP}.
\end{align}
The recurrence relation of $P(n,k)$ in Equation \eqref{eq:recursiveP} looks similar to that of the Eulerian number $\eulerian{n}{k}$ from Equation \eqref{eq:eulrec}. 

If we assume 
\[
P(n-1,k) = \frac{\eulerian{n-1}{k}}{(n-1)!},
\] 
for $1 \leq k \leq n-2$, then,
\begin{align*}
P(n,k) &= q_{n,k}P(n-1,k-1) + p_{n,k}P(n-1,k) \\
 &= \frac{n-k}{n} \cdot \frac{\eulerian{n - 1}{k - 1}}{(n-1)!} + \frac{k+1}{n} \cdot \frac{\eulerian{n - 1}{k}}{(n-1)!},\\
  &=\frac{(n - k)\eulerian{n - 1}{k - 1} + (k + 1)\eulerian{n - 1}{k}}{n!},\\
  &= \frac{\eulerian{n}{k}}{n!},
\end{align*}
where the final equality follows from the recurrence for Eulerian numbers in Equation \eqref{eq:recursiveP}.
This proves our first theorem.

\begin{thm}\label{thm:eulprob}
For any $n \geq 1$ and $0 \leq k \leq n - 1$,
\[
P(n,k) = \frac{\eulerian{n}{k}}{n!}.
\]
\end{thm} 

Our first main result is that $P(n,k)$, the probability of arriving at state $s(n,k)$, is modeled by the Eulerian probability distribution. In other words, we found a nice result for the probability that after emptying a queue of $n$ particles, $k$ of them settle to the right of the origin. With question (1) answered, we now move onto question (2): the expected amount of time (i.e., the number of coin tosses) it takes for $n$ particles to settle. 

\section{Expected run time}

We next turn to the expected run time of our coin-tossing game, i.e., we want to know how many tosses it will take for $n$ particles to settle. Recall from the introduction that $E_n$ denotes the expected number of tosses needed to reach a state with $n$ occupied sites. Notice that $E_1 = 0$, and $E_2 = 1$ because flipping heads or tails will cause a site to become occupied immediately. It gets more complicated once we start looking at cases for $n \geq 3$. In the introduction, we gave an estimate for $E_n$ based on our Monte Carlo simulation. We now turn this into our next theorem.

\begin{thm}\label{thm:tosstosn}
For any $n \geq 1$, the expected number of tosses needed to reach a state with $n$ occupied sites is 
\begin{equation}\label{eq:actualtoss}
E_n = \frac{1}{12} n^3 + \frac{1}{12} n^2.
\end{equation}
\end{thm}

If we compute the first differences of this sequence ($E_{n} - E_{n-1}$), we find the change in expected tosses. Let $\Delta E_n$ denote the expected number of tosses it takes for the $n^{th}$ particle to perform a random walk and settle, given $n-1$ sites are already occupied.  
Assuming Equation \eqref{eq:actualtoss}, we find
\[
\Delta E_n = \frac{1}{12} n^3 + \frac{1}{12} n^2 - \left(\frac{1}{12} (n-1)^3 + \frac{1}{12} (n-1)^2\right) = \frac{1}{4}n^2 - \frac{1}{12}n.
\]
Thus to obtain our run time result in Theorem \ref{thm:tosstosn}, it suffices to show that the number of expected tosses for the $n^{th}$ particle, once dropped, to settle in an unoccupied site is
\begin{equation}\label{eq:newtossquad}
\Delta E_n = \frac{1}{4}n^2 - \frac{1}{12}n.
\end{equation}

But by definition, the expected number of new steps is 
\[
\Delta E_n = \sum_{k=1}^{n-1}P(n-1,k-1) \cdot E(n-1,k-1),
\]
where $E(n-1,k-1)$ denotes the expected number of tosses to go from state $s(n-1,k-1)$ to $s(n,k)$ or $s(n,k-1)$. By Theorem \ref{thm:eulprob}, we can write
\[
 \Delta E_n = \sum_{k=1}^{n-1} \frac{\eulerian{n-1}{k-1}}{(n-1)!} \cdot E(n-1,k-1) = \frac{1}{(n-1)!} \sum_{k=1}^{n-1} \eulerian{n-1}{k-1} E(n-1,k-1),
\]
i.e., $\Delta E_n$ is a weighted sum of Eulerian numbers. We will show $E(n-1,k-1)$ is simply $k(n-k)$, which will help us to evaluate the sum and verify \eqref{eq:newtossquad}, and hence Theorem \ref{thm:tosstosn}.

\subsection{Escape time for a random walk}

Consider that the game is in state $s(n-1,k-1)$ pictured below, 
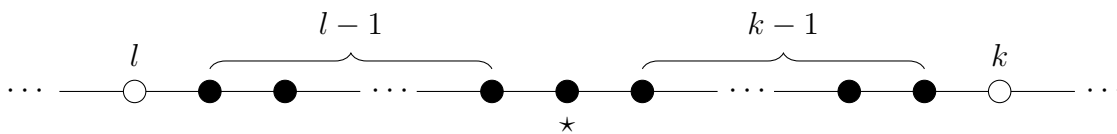
\[
\begin{tikzpicture}
\tikzstyle{wn}=[inner sep =3, fill=white, draw=black, circle];
\tikzstyle{bn}=[inner sep =3, fill=black, draw=black, circle];
\draw (0,0)--(1,0) node[wn] {} --(2,0) node[bn] {}--(3,0) node[bn]{} -- (4,0);
\draw (4,0) node[right] {$\cdots$};
\draw (4.75,0)--(5.75,0) node[bn] {} --(6.75,0) node[bn] {} -- (7.75,0) node[bn] {} -- (8.75,0); 
\draw (8.75,0) node[right]{$\cdots$};
\draw (9.5,0)--(10.5,0) node[bn] {}--(11.5,0) node[bn]{} -- (12.5,0) node[wn]{} -- (13.5,0);
\draw (0,0) node[left] {$\cdots$};
\draw (13.5,0) node[right] {$\cdots$};
\draw (6.75,-.2) node[below] {$\star$}; \draw (12.5,.2) node[above] {$k$}; \draw (1,.2) node[above] {$l$};
\draw [decoration={brace,amplitude=6pt},decorate] (7.75,0.3) -- (11.5,0.3) node[midway,above=8pt] {$k-1$};
\draw [decoration={brace,amplitude=6pt},decorate] (2,0.3) -- (5.75,0.3) node[midway,above=8pt] {$l-1$};
\end{tikzpicture}
\]
where $k+l = n$. We would like to know how to compute $E(n-1,k-1),$ i.e., how many coin tosses it will take for the game to reach either state $s(n,k)$ (activating on the right) or state $s(n,k-1)$ (activating on the left). To help us, recall that for the simple random walk, the expected time to reach the boundary of an interval $[-b,a]$ is $ab$. See \cite{Sauer}. We translate this escape time result into the language of our game with the following claim. 
\begin{cl}\label{cl:esctime}
For any $k,l \geq 1$ such that $k+l=n$, the expected number of tosses to go from state $s(n-1,k-1)$ to $s(n,k)$ or $s(n,k-1)$ is $E(n-1,k-1) = kl = k(n-k)$.
\end{cl}

For example, consider state $s(5,3)$ below.
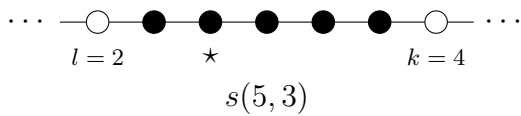
\[
\begin{tikzpicture}
\tikzstyle{wn}=[inner sep =3, fill=white, draw=black, circle];
\tikzstyle{bn}=[inner sep =3, fill=black, draw=black, circle];

\draw (9,-2)--(9.5,-2) node[wn] {} -- (10.25,-2) node[bn] {} -- (11,-2) node[bn] {} -- (11.75,-2) node[bn] {} -- (12.5,-2) node[bn]{} -- (13.25,-2) node[bn] {} -- (14,-2) node[wn] {} -- (14.5,-2); 
\draw (9,-2) node[left] {$\cdots$}; \draw (14.5,-2) node[right] {$\cdots$};
\draw (11,-2.2) node[below] {$\star$}; \draw (9.5,-2.2) node[below] {\scriptsize $l=2$}; \draw (14,-2.2) node[below] {\scriptsize $k=4$};
\draw (11.75,-2.65) node[below] {$s(5,3)$};
\end{tikzpicture}
\]

As shown, $l=2$ and $k=4$ in $s(5,3)$. We can calculate, using Claim \ref{cl:esctime}, the expected number of tosses it takes to go from each state with $n=5$ to a state with $n=6$ occupied sites. Consider the chart below:
\[
\begin{centering}
\begin{tabular}{c | c | c | c} 
State & $k$ & $l$ & $E(5,k-1) = kl$ \\
\hline
$s(5,0)$ & $1$ & $5$ & 5 \\
$s(5,1)$ & $2$ & $4$ & 8 \\
$s(5,2)$ & $3$ & $3$ & 9 \\
$s(5,3)$ & $4$ & $2$ & 8 \\
$s(5,4)$ & $5$ & $1$ & 5 \\
\end{tabular}
\end{centering}
\]
If we know there are five occupied sites (but not which 5), and we want to know how many tosses it would take for the sixth particle to settle, we calculate a weighted sum using our Eulerian probabilities. This is precisely our definition for $\Delta E_6$. 
\[
\Delta E_6 = \left(\frac{1}{120}\right)(5) + \left(\frac{26}{120}\right)(8) + \left(\frac{66}{120}\right)(9) + \left(\frac{26}{120}\right)(8) + \left(\frac{1}{120}\right)(5) = \frac{17}{2}. 
\]
This agrees with Equation \eqref{eq:newtossquad}, 
\[ 
\frac{1}{4} \left(6^2\right) - \frac{1}{12} \left(6\right) = \frac{17}{2}.
\]

In general, since $P(n,k) = \frac{\eulerian{n}{k}}{n!}$ by Theorem \ref{thm:eulprob} and $E(n-1,k-1) = kl$ by Claim \ref{cl:esctime}, we have the following Corollary of Theorem \ref{thm:eulprob} and Claim \ref{cl:esctime}.

\begin{cor}\label{cor:deltaen}
For any $n \geq 3$, the expected number of tosses for the $n^{th}$ particle to settle, given $n-1$ sites are occupied is 
\begin{align*}
\Delta E_n &= \sum_{k=1}^{n-1} \frac{\eulerian{n-1}{k-1}}{(n-1)!} \cdot k(n-k) \\
&= \frac{1}{(n-1)!} \sum_{k=1}^{n-1} \eulerian{n-1}{k-1} \cdot k(n-k).
\end{align*}
\end{cor}

It remains to see that this also equals the formula in Equation \eqref{eq:newtossquad}. We will do so using generating function techniques.

\subsection{Eulerian polynomials and generating functions}

Throughout this section we use the technique of \emph{generating functions} which is explored in \cite{Wilf}, and in particular \cite{Petersen} for Eulerian numbers. The basic idea is to study power series whose coefficients are the numbers we are interested in, and then to use properties of these functions to deduce properties of the coefficients. 

To this end, we define the $n^{th}$ bivariate Eulerian polynomial, denoted $A_{n}(s,t)$, as 
\[
A_{n}(s,t) = \sum_{k=1}^{n} \eulerian{n}{k-1} s^{k} t^{n+1-k}, 
\]
e.g., $A_5 (s,t) = st^5 + 26s^2t^4 + 66s^3t^3 + 11s^4t^2 + s^5t$. 

We can get $\Delta E_n$ from $A_{n-1} (s,t)$ as follows. First, we differentiate with respect to $s$ and $t$,
\begin{align*}
\frac{\partial}{\partial s} \frac{\partial}{\partial t} \left[A_{n-1}(s,t)\right] &= \frac{\partial}{\partial s} \frac{\partial}{\partial t} \left[\sum_{k=1}^{n-1} \eulerian{n-1}{k-1}s^kt^{n-k} \right], \\
&= \sum_{k=1}^{n-1} \eulerian{n-1}{k-1}k(n-k)s^{k-1}t^{n-k-1},
\end{align*}
and when $s=1$ and $t=1$, this is 
\begin{equation}\label{eq:wgtsum}
\sum_{k=1}^{n-1} \eulerian{n-1}{k-1}k(n-k) = (n-1)! \cdot \Delta E_n.
\end{equation}

For example, consider the case for $n=6$. We differentiate $A_5 (s,t)$ with respect to $s$ and $t$ and then set $s=t=1$,
\[
\frac{\partial}{\partial s} \frac{\partial}{\partial t} \left[A_{5}(s,t)\right]_{s=t=1} = (1\cdot 5) + 26(2 \cdot 4) + 66(3\cdot 3) + 26(4 \cdot 2) + (5\cdot 1).
\]
Notice that this is precisely the weighted sum for $n=6$ from Equation \eqref{eq:wgtsum},
\[
\sum_{k=0}^{4} \eulerian{5}{k-1} k(6-k) = 5! \cdot \Delta E_6.
\]

Before we can use this idea to prove Equation \eqref{eq:newtossquad} holds in general, we consider the exponential generating function for the Eulerian polynomial $A_n(s,t)$, denoted $F(s,t,z)$. By definition, we have   
\begin{equation}\label{eq:expgenf}
F(s,t,z) = \sum_{n \geq 1}A_{n}(s,t) \frac{z^n}{n!}.
\end{equation}
That is, expanding in $z$, we find 
\[
F(s,t,z) =  stz + (s^2t + st^2 ) \frac{z^2}{2} + (s^3t + 4s^2t^2 + st^3) \frac{z^3}{3!} + \cdots .
\]

To find a useful formula for $F(s,t,z),$ we work to the bivariate case from the univariate case, which is better known. 
Notice that 
\begin{equation}\label{eq:uni2bi}
A_n (s,t) = \sum_{k=1}^n \eulerian{n}{k-1}s^kt^{n+1-k} = t^{n+1}\sum_{k=1}^n \eulerian{n}{k-1} \left(\frac{s}{t}\right)^k= t^{n+1} A_n \left(\frac{s}{t}\right),
\end{equation}
where $A_n(x) = \sum_{k=1}^{n} \eulerian{n}{k-1} x^k$ is the single-variable Eulerian polynomial.

The exponential generating function for $A_n(x)$, denoted $G(x,z)$, is 
\[
G(x,z) = \sum_{n\geq 1} A_n (x) \frac{z^n}{n!} = \frac{x\left(1-e^{z(x-1)}\right)}{e^{z(x-1)} -x}.
\]
See \cite{Petersen}. Thus, we see that 
\begin{align}
F(s,t,z) &= \sum_{n \geq 1} A_n (s,t) \frac{z^n}{n!} \nonumber \\
&= t \cdot \sum_{n \geq 1} A_n \left(\frac{s}{t}\right) \frac{(tz)^n}{n!} \nonumber \\
&= t \cdot G\left(\frac{s}{t},tz\right) \nonumber \\
&= \frac{st\left(1- e^{z(s-t)}\right)}{te^{z(s-t)} -s} \label{eq:fstz}. 
\end{align}

With the function $F(s,t,z)$ in Equation \eqref{eq:fstz}, we can calculate $\Delta E_n$ simultaneously for all $n$ through the same steps of differentiation and substitution, i.e., 
\begin{align*}
\frac{\partial}{\partial s} \frac{\partial}{\partial t} [F(s,t,z)]_{s=t=1} &= \sum_{n \geq 1} \left(\frac{\partial}{\partial s} \frac{\partial}{\partial t} [A_n(s,t)]_{s=t=1}\right) \frac{z^n}{n!}, \\
&= \sum_{n \geq 1} n!\cdot \Delta E_{n+1} \frac{z^n}{n!} = \sum_{n \geq 1} \Delta E_{n+1} z^n.
\end{align*}
Taking the deriviative at $s=t=1$ of the formula for $F(s,t,z)$ given in \eqref{eq:fstz}, we find
\[
\frac{z^4 - 4z^3 + 6z^2 -6z}{6(z-1)^3} = \sum_{n \geq 1} \Delta E_{n+1} \cdot z^n. 
\]

Note that from Equation \eqref{eq:newtossquad}, we would like to see
\begin{equation}\label{eq:en1}
\Delta E_{n+1} = \frac{1}{4} \left(n+1\right)^2 - \frac{1}{12} \left(n+1\right) = \frac{1}{4} n^2 + \frac{5}{12} n + \frac{1}{6}. 
\end{equation}

We will now use generating function techniques to show that 
\[
\frac{z^4 - 4z^3 + 6z^2 -6z}{6(z-1)^3} = \sum_{n\geq 1} \left(\frac{1}{4} n^2 + \frac{5}{12} n + \frac{1}{6}\right) z^n.
\]
 
By differentiating the geometric series $1/(1-z) = 1+z+z^2 + \cdots$ twice (and dividing by $2$), we know that 
\begin{equation}\label{eq:genfunc}
\frac{1}{(1-z)^3} = \sum_{n \geq 0} \binom{n+2}{2} z^n = \binom{2}{2} + \binom{3}{2}z + \binom{4}{2} z^2 + \cdots.
\end{equation}

We split our generating function for expectation into four pieces:
\[
\frac{z^4 - 4z^3 + 6z^2 -6z}{6(z-1)^3} = \frac{6z}{6(1-z)^3} - \frac{6z^2}{6(1-z)^3} + \frac{4z^3}{6(1-z)^3} - \frac{z^4}{6(1-z)^3},
\]
and using the formula from Equation \eqref{eq:genfunc}, we find 
\begin{align*}
\frac{6z}{6(1-z)^3} &= z \cdot \sum_{n \geq 0} \binom{n+2}{2} z^n = \sum_{n\geq 1} \binom{n+1}{2} z^n, \\
-\frac{6z^2}{6(1-z)^3} &=- z^2 \cdot \sum_{n \geq 0} \binom{n+2}{2} z^n = -\sum_{n\geq 2} \binom{n}{2} z^n, \\
\frac{4z^3}{6(1-z)^3} &=\frac{2z^3}{3} \cdot \sum_{n \geq 0} \binom{n+2}{2} z^n = \sum_{n\geq 3} \frac{2}{3} \binom{n-1}{2}z^n,\\
- \frac{z^4}{6(1-z)^3} &= - \frac{z^4}{6} \cdot \sum_{n \geq 0} \binom{n+2}{2} z^n = \sum_{n \geq 4} \frac{1}{6} \binom{n-2}{2} z^n.
\end{align*}

Summing these, we get
\begin{align*}
 \frac{z^4 - 4z^3 + 6z^2 -6z}{6(z-1)^3} &= z + 2z^2 + \frac{11}{3}z^3 \\
  & \quad + \sum_{n\geq 4}\left[\binom{n+1}{2} -\binom{n}{2} + \frac{2}{3}\binom{n-1}{2} + \frac{1}{6}\binom{n-2}{2} \right] z^n.
\end{align*}

Using the binomial formula to simplify the general term for $n\geq 4$, we find
\[
\binom{n+1}{2} - \binom{n}{2} + \frac{2}{3}\binom{n-1}{2} - \frac{1}{6}\binom{n-1}{2} = \frac{1}{4}n^2 + \frac{5}{12}n + \frac{1}{6},
\]
which establishes Equation \eqref{eq:en1}, and thus Theorem \ref{thm:tosstosn}.

\section{Further remarks}

We also have Monte Carlo results about question (3), the probability that $n$ sites are occupied given a coin is tossed $N$ times (with a large queue of particles). From the data,
it appears that the probabilities have a common denominator of $2^{n-1}$. After clearing denominators and rounding the values to the nearest integer, the results for $1\leq n \leq 7$ are as follows:
\[
\begin{array}{c | c c c c c}
N \backslash n & 2 & 3 & 4 & 5 \\
\hline
1 & 1 \\
2 & 1 & 1 \\
3 & 1 & 3 \\
4 & 1 & 4 & 3 \\
5 & 1 & 9 & 6 \\
6 & 1 & 10 & 18 & 3 \\
7 & 1 & 21 & 32 & 10 
\end{array}
\]
We did not find this array in the On-Line Encyclopedia of Integer Sequences \cite{Sloane}, though the first two columns have simple explanations. Also, we remark that the first nonzero entry in column $n$ occurs when $N = \lceil \frac{n}{2} \rceil \cdot \lfloor \frac{n}{2} \rfloor$.
 
We finish with some thoughts and observations about question (4) about what happens when we no longer have a simple random walk. Using linear algebra, one can analyze the one-dimensional random walk for any fixed choice of probabilities of a particle moving left or right at a given site. The coin toss model we used is just the case where the probability is $\frac{1}{2}$ at each site. We have experimented with other choices of probabilities, but without interesting results, except in the case of a biased coin toss. Here each site has the same fixed probability, but that probability is not necessarily $\frac{1}{2}$. In this case we can use analysis similar to that in this paper to obtain similar results, which we summarize now. 

Let $p$ be the probability of a particle going to the right and $q$ the probability of going to the left, and let $\rho=p/q$. The biased Gambler's ruin gives transition probabilities of $[k]/[n]$ where
\[
 [n] = \frac{1-\rho^n}{1-\rho} = 1+\rho+\rho^2 + \cdots + \rho^{n-1}.
\]
Then the same arguments we used in the fair coin case will show that 
\begin{equation}
 P(n,k) = \frac{ \eulerian{n}{k}^{\maj}}{[n]!},
\end{equation}
where $[n]! = [n][n-1]\cdots[2][1]$, and $\eulerian{n}{k}^{\maj}$ counts permutations of length $n$ with $k$ descents, weighted according to a power of $\rho$ given by a statistic called \emph{major index}. That is, 
\[
\eulerian{n}{k}^{\maj} = \sum_{\substack{w \in S_n, \\ \des(w) =k }} \rho^{\maj(w)}.
\]
This gives a refinement of the Eulerian distribution in terms of $\rho$, as explained in Chapter 6 of \cite{Petersen}. Studying questions like the limiting distribution and run time when $\rho\neq 1$ is the starting point for future work on this problem.

\section*{Acknowledgements}

The author is grateful for support from the Undergraduate Summer Research Program from the College of Science and Health at DePaul University, and the guidance of her advisor Kyle Petersen. Thanks also to Lionel Levine for suggesting the problem, and thanks to Jim Propp for comments on an earlier draft.

\end{document}